\documentclass[12pt]{article}
\usepackage{amssymb}

\newcommand{\real}{{\bf R}}
\newcommand{\complex}{{\bf C}}
\newcommand{\intplus}{{\bf N}}
\newcommand{\allint}{{\bf Z}}

\renewcommand{\d}{\,{\rm d}}            
\newcommand{\D}{{\rm d}}                
\renewcommand{\Re}{{\rm Re\,}}

\renewcommand{\epsilon}{\varepsilon}
\renewcommand{\phi}{\varphi}

\newcommand{\cL}{{\cal L}}
\newcommand{\cN}{{\cal N}}
\newcommand{\cM}{{\cal M}}
\newcommand{\cO}{{\cal O}}

\newcommand{\cS}{{\cal S}}

\newcommand{\ee}{{\bf e}}
\newcommand{\uu}{{\bf u}}
\newcommand{\vv}{{\bf v}}

\newcommand{\oone}{{\bf 1}}

\newcommand{\UU}{{\bf U}}
\newcommand{\OO}{{\bf \Omega}}

\newtheorem{theorem}{Theorem}[section]
\newtheorem{lemma}[theorem]{Lemma}

\newtheorem{proposition}[theorem]{Proposition}
\newtheorem{corollary}[theorem]{Corollary}
\newtheorem{remark}[theorem]{Remark}

\newcommand{\reff}[1]{(\ref{#1})}

\newcommand{\inttwo}{\int_{\real^2}}
\newcommand{\proof}{{\noindent \bf Proof:\ }}
\newcommand{\half}{\textstyle{\frac{1}{2}}}

\def\build#1_#2^#3{\mathrel{
  \mathop{\kern 0pt#1}\limits_{#2}^{#3}}}
\def\QED{\mbox{}\hfill$\Box$}

\usepackage{amsfonts}
\begin{document}

\title{Existence and stability of asymmetric Burgers vortices}

\author{Thierry Gallay \\ Institut Fourier \\
Universit\'e de Grenoble I \\ BP 74 \\
38402 Saint-Martin-d'H\`eres \\France
\and
C. Eugene Wayne \\
Department of Mathematics \\
and  Center for BioDynamics \\
Boston University \\
111 Cummington St.\\
Boston, MA 02215, USA}

\date{\normalsize March 15, 2005}

\maketitle
\begin{abstract}
Burgers vortices are stationary solutions of the three-dimensional 
Navier-Stokes equations in the presence of a background straining
flow. These solutions are given by explicit formulas only when the 
strain is axisymmetric. In this paper we consider a weakly asymmetric 
strain and prove in that case that non-axisymmetric vortices 
exist for all values of the Reynolds number. In the limit 
of large Reynolds numbers, we recover the asymptotic results of 
Moffatt, Kida \& Ohkitani \cite{moffatt:1994}. We also show that 
the asymmetric vortices are stable with respect to localized
two-dimensional perturbations. 
\end{abstract}

\section{Introduction}
\label{intro} 

Localized structures such as vortex sheets or tubes play a prominent
role in the dissipation of energy in three-dimensional turbulent
flows. It is believed that these dissipative structures take place due
to the interplay of two basic mechanisms: amplification of vorticity
due to stretching, and diffusion through the action of viscosity
\cite{taylor:1938}. A typical example that exhibits both features is
the familiar Burgers vortex \cite{burgers:1948}, an explicit solution
of the three-dimensional Navier-Stokes equations in the presence of an
axisymmetric background straining flow. In real flows, however, the
local strain has no reason of being axisymmetric, and as a matter of
fact the vortex tubes observed in numerical simulations usually
exhibit a truly elliptical core region. It is therefore important to
study the analogue of the Burgers vortex when the straining flow is 
asymmetric, although no explicit expression is available in that case. 

Using a double series expansion, Robinson and Saffman
\cite{robinson:1984} formally established the existence of an
asymmetric vortex for small values of the Reynolds number $R$ and of
the asymmetry parameter $\lambda$. This solution was also studied
numerically for larger $\lambda$ (up to 3/4) and $R$ (up to 100). On
the other hand, an asymptotic expansion for large Reynolds numbers was
performed by Moffatt, Kida and Ohkitani \cite{moffatt:1994}, see also
\cite{jimenez:1996}. Their results indicate that an equilibrium
stretched vortex should exist for all values of $\lambda \in (0,1)$
and $R > 0$ such that $\lambda/R \ll 1$. Interesting features of these
solutions, such as the shape of isovorticity contours and the spatial
distribution of energy dissipation, were also studied in detail.
Finally, the stability of symmetric or non-symmetric vortices is an
important issue which has attracted a lot of attention in recent
years. Roughly speaking, the stability with respect to two-dimensional
perturbations (i.e., perturbations which are independent of the axial
coordinate) is well understood \cite{robinson:1984,prochazka:1995,
prochazka:1998, gallay:2004}, but only partial results have been 
obtained in the general case where arbitrary three-dimensional 
perturbations are allowed \cite{rossi:1997,crowdy:1998,eloy:1999,
schmid:2004, gallay:2005}.

In this paper, we prove (rigorously) that non-axisymmetric Burgers
vortices exist for all values of the Reynolds number, provided the
asymmetry parameter is sufficiently small. In particular, taking the
limit $R \to \infty$, we recover exactly the asymptotic results of
Moffatt, Kida and Ohkitani \cite{moffatt:1994}. We also show that
these vortices are stable with respect to spatially localized
two-dimensional perturbations. Existence for larger values of the
asymmetry parameter and stability with respect to three-dimensional
perturbations are difficult questions, which have been solved so far
for small Reynolds numbers only \cite{gallay:2005}.

We now describe our results in more detail. We consider an incompressible
viscous fluid filling the whole space $\real^3$, and we suppose that
the velocity field is a two-dimensional perturbation of a linear
straining flow, namely
$$
  \UU(x_1,x_2,x_3,t) \,=\, \pmatrix{\gamma_1 x_1 \cr \gamma_2 x_2 \cr
  \gamma_3 x_3} + \pmatrix{u_1(x_1,x_2,t)\cr u_2(x_1,x_2,t) \cr 0}~,
$$
where $\gamma_1,\gamma_2,\gamma_3$ are reals constants satisfying 
$\gamma_1 + \gamma_2 + \gamma_3 = 0$. Throughout this paper we assume
that
\begin{equation}\label{gamdef}
  \gamma_1 \,=\, -\frac{\gamma}2 (1+\lambda)~, \quad 
  \gamma_2 \,=\, -\frac{\gamma}2 (1-\lambda)~, \quad 
  \gamma_3 \,=\, \gamma~, 
\end{equation}
for some $\gamma > 0$ and some $\lambda \in [0,1)$. Thus $\gamma_3$
is the only positive principal rate of strain, and the straining flow
is axisymmetric if and only if $\lambda = 0$. The case of a biaxial 
strain ($\lambda > 1$), which is also important for applications
in turbulence, will not be considered here. 

The vorticity $\OO = \nabla \times \UU$ is aligned with the vertical
axis and depends only on the horizontal variable, namely
$$
  \OO(x_1,x_2,x_3,t) \,=\, \pmatrix{0\cr 0 \cr \omega(x_1,x_2,t)}~,
  \quad \hbox{where} \quad \omega \,=\, \partial_1 u_2 - \partial_2 u_1~.
$$
The evolution equation for $\omega$ reads
\begin{equation}\label{omeq1}
  \partial_t \omega + (u_1 + \gamma_1 x_1)\partial_1 \omega + 
  (u_2 + \gamma_2 x_2)\partial_2 \omega \,=\, \nu \Delta \omega + 
  \gamma_3 \omega~,
\end{equation}
where $\nu > 0$ is the kinematic viscosity of the fluid. Since 
$\partial_1 u_1 + \partial_2 u_2 = 0$ and $\partial_1 u_2 - 
\partial_2 u_1 = \omega$, the rotational velocity $\uu = (u_1,u_2)$
can be recovered from $\omega$ via the two-dimensional Biot-Savart law
\begin{equation}\label{BS}
  \uu(x,t) \,=\, \frac{1}{2\pi} \inttwo \frac{(x-y)^\perp}{|x-y|^2}
  \omega(y,t)\d y~, \quad x = (x_1,x_2) \in \real^2~,
\end{equation}
where $x^\perp = (-x_2,x_1)$ and $|x|^2 = x_1^2 + x_2^2$. 

In the axisymmetric case $\lambda = 0$, equation \reff{omeq1} 
has a family of explicit time-independent solutions:
\begin{equation}\label{burgdef}
  \omega^\Gamma \,=\, \frac{\Gamma}{\delta^2}\,G\Bigl(\frac{x}{\delta}
  \Bigr)~, \quad \uu^\Gamma \,=\, \frac{\Gamma}{\delta}\,\vv^G
  \Bigl(\frac{x}{\delta}\Bigr)~,
\end{equation}
where $\Gamma \in \real$, $\delta = (\nu/\gamma)^{1/2}$, and
\begin{equation}\label{GvG}
  G(x) \,=\, \frac{1}{4\pi} \,e^{-|x|^2/4}~,\quad
  \vv^G(x) \,=\, \frac{1}{2\pi}\frac{x^\perp}{|x|^2}
  \Bigl(1 -  e^{-|x|^2/4}\Bigr)~, \quad x \in \real^2\ .
\end{equation}
These are the well-known (axisymmetric) Burgers vortices. The family
is indexed by the parameter $\Gamma = \inttwo \omega^\Gamma\d x$, 
which represents the circulation of $\uu^\Gamma$ at infinity. The
Reynolds number associated to the Burgers vortex with circulation
$\Gamma$ can be defined \cite{moffatt:1994} as
$$
   R \,=\, \frac{|\Gamma|}{\nu}~.
$$

The aim of this paper is to study the analogue of the Burgers vortices 
when the straining flow is not axisymmetric. The expressions of these
asymmetric vortices will be greatly simplified if we use the natural
lengthscale $\delta = (\nu/\gamma)^{1/2}$ and timescale $\tau = 
1/\gamma$ defined by the viscosity and the strain. We thus replace the 
variables $x,t$ and the functions $\uu,\omega$ with the dimensionless 
quantities
$$
  \tilde x \,=\, \frac{x}{\delta}~, \quad
  \tilde t \,=\, \frac{t}{\tau}~, \quad 
  \tilde \uu \,=\, \frac{\tau \uu}{\delta}~, \quad 
  \tilde \omega \,=\, \tau\omega~.
$$
Dropping the tildes for convenience, we see that the new functions
$\omega,\uu$ satisfy \reff{omeq1} with $\gamma = \nu = 1$, namely
\begin{equation}\label{omeq2}
  \partial_t \omega + \uu \cdot \nabla \omega \,=\, \cL \omega
  + \lambda \cM \omega~, 
\end{equation}
where
\begin{equation}\label{LMdef}
  \cL \,=\, \Delta + \half x\cdot\nabla + 1~, \quad
  \cM \,=\, \half(x_1\partial_1 - x_2\partial_2)~.
\end{equation}
It is easily verified that $\cL G = 0$ and $\vv^G \cdot \nabla G = 0$, 
hence in the symmetric case $\lambda = 0$ the Burgers vortex 
$\omega^\alpha = \alpha G$, $\uu^\alpha = \alpha\vv^G$ is indeed 
a stationary solution of \reff{omeq2} for any $\alpha \in \real$. 
The Reynolds number associated to this flow is simply $R = |\alpha|$.

To formulate our results, we introduce appropriate function spaces.
Let $X$ be the (real) Hilbert space
\begin{equation}\label{Xdef}
  X \,=\, \Bigl\{w \in L^2(\real^2) \,\Big|\, G^{-1/2}w \in L^2(\real^2)\,,
  ~\inttwo w\d x = 0\Bigr\}~,
\end{equation}
equipped with the scalar product
$$
  (w_1,w_2)_X \,=\, \inttwo G(x)^{-1} w_1(x)w_2(x)\d x~.
$$
We also define the subspace $Y = \{w \in X \,|\, \partial_i w \in X
\hbox{ for }i = 1,2\}$ equipped with the natural scalar product
$$
  (w_1,w_2)_Y \,=\, \inttwo G(x)^{-1}\Bigl(w_1(x)w_2(x) + 
  \nabla w_1(x)\cdot\nabla w_2(x)\Bigr)\d x~.
$$
Except for the zero mean condition, the space $X$ is just a weighted
$L^2$ space (with Gaussian weight) and $Y$ is the corresponding 
Sobolev space. We can now state our first result:

\begin{theorem}\label{thm1}
There exist $\lambda_0 > 0$ and $K_0 > 0$ such that, for all
$\lambda \in [0,\lambda_0]$ and all $\alpha \in \real$, equation
\reff{omeq2} has a unique stationary solution
$\omega^{\alpha,\lambda},\uu^{\alpha,\lambda}$ such that
$\|\omega^{\alpha,\lambda} - \alpha G\|_Y \le K_0$. Moreover, 
$\|\omega^{\alpha,\lambda} - \alpha G\|_Y \le K_0(\lambda/\lambda_0)
|\alpha|/(1{+}|\alpha|)$.
\end{theorem}

This theorem shows that, if the asymmetry parameter $\lambda$ is 
sufficiently small, equation~\reff{omeq2} has a family of equilibria
$\omega^{\alpha,\lambda}$ indexed by the circulation number 
$\alpha = \inttwo \omega^{\alpha,\lambda}\d x$. The solution
$\omega^{\alpha,\lambda}$ is locally unique, has a Gaussian decay 
at infinity, and converges to the Burgers vortex $\alpha G$ as 
$\lambda \to 0$. Further properties of these asymmetric vortices 
will be established in Section~\ref{existence}. For instance, 
$\omega^{\alpha,\lambda}(x)$ is a smooth function of $x \in \real^2$, 
$\alpha \in \real$, and $\lambda \in [0,\lambda_0]$. Moreover, 
$\omega^{\alpha,\lambda} > 0$ if $\alpha > 0$,
$\omega^{\alpha,\lambda} < 0$ if $\alpha < 0$, and 
$\omega^{\alpha,\lambda} \equiv 0$ if $\alpha = 0$.

Theorem~\ref{thm1} will be proved by a classical perturbation
argument. The only remarkable point is that this argument can be
applied uniformly for all $\alpha \in \real$. In particular, for fixed
$\lambda$, we can investigate the limit of large Reynolds numbers $R =
|\alpha| \to \infty$. The last inequality in the theorem asserts that
$\alpha^{-1}\omega^{\alpha,\lambda} = G + \cO(\lambda/R)$ as $R \to
\infty$, hence we recover the observation by Moffatt, Kida and
Ohkitani \cite{moffatt:1994} that the asymptotic profile is always the
Gaussian $G$, even in the asymmetric case $\lambda > 0$. (Note that,
unlike in \cite{moffatt:1994}, we do not need to {\em assume} that
the asymptotic profile is radially symmetric.) Moreover, the
deviation from the limiting profile is proportional to $\lambda/R$ at
leading order, as established in \cite{moffatt:1994}. A rigorous
expansion up to second order in $\lambda$ and $R^{-1}$ will be
performed in Section~\ref{existence}, see Eq.\reff{largeR} below. 

Our second result shows that the asymmetric Burgers vortex
is asymptotically stable stable with respect to perturbations
in $X$. As in Theorem~\ref{thm1}, this property holds uniformly
for all $\alpha \in \real$ and $\lambda \in [0,\lambda_1]$, 
for some $\lambda_1 > 0$ (possibly smaller than $\lambda_0$). 
Remark that there is no loss of generality in assuming that
the perturbations have zero mean, because if $\int \tilde 
\omega \d x = \beta \neq 0$ then $\omega^{\alpha,\lambda} +
\tilde\omega$ is a zero mean perturbation of the (modified) vortex 
$\omega^{\alpha+\beta,\lambda}$. 

\begin{theorem}\label{thm2}
Given any $\mu \in (0,1/2)$, there exist $\lambda_1 > 0$ and 
$\epsilon > 0$ such that, for all $\lambda \in [0,\lambda_1]$ and 
all $\alpha \in \real$, the following holds. For all initial data 
$\omega_0$ with $\|\omega_0 - \omega^{\alpha,\lambda}\|_X \le \epsilon$, 
equation \reff{omeq2} has a unique global solution $\omega(x,t)$
such that $\omega - \omega^{\alpha,\lambda} \in C^0([0,+\infty),X)$.
Moreover, 
\begin{equation}\label{asymptotics}
  \|\omega(\cdot,t) - \omega^{\alpha,\lambda}\|_X \,\le\, 
  \|\omega_0 - \omega^{\alpha,\lambda}\|_X\,e^{-\mu t}~, \quad
  \hbox{for all } t\ge 0~. 
\end{equation}
\end{theorem}

The proof of Theorem~\ref{thm1} is based on the following ideas.
Given $\lambda \in (0,1)$ and $\alpha \in \real$, we look for 
stationary solutions of \reff{omeq2} of the form $\omega = \alpha G 
+ w$, $\uu = \alpha \vv^G + \vv$, where $\inttwo w\d x = 0$ and 
$\vv$ is the velocity field obtained from $w$ via the Biot-Savart
law \reff{BS}. The equation for $w$ reads:
\begin{equation}\label{weq1}  
  \alpha(\vv^G\cdot\nabla w + \vv\cdot\nabla G) + \vv\cdot\nabla w
  \,=\, \cL w + \lambda\cM (\alpha G + w)~.
\end{equation}
Let $\Lambda$ be the integro-differential operator defined by
\begin{equation}\label{Lambdadef}
  \Lambda w \,=\, \vv^G \cdot\nabla w + \vv\cdot\nabla G~.
\end{equation}
It is shown in \cite{gallay:2004} that, for any $\alpha \in \real$, 
the spectrum of $\cL - \alpha\Lambda$ acting on $X$ is contained
in the half-plane $\{z \in \complex \,|\, \Re(z)\le -\half\}$. 
In particular $\cL -\alpha\Lambda$ is invertible, and \reff{weq1}
can be rewritten as
\begin{equation}\label{weq2}  
  w \,=\, (\cL - \alpha\Lambda)^{-1}\Bigl(\vv\cdot\nabla w - 
 \lambda\cM (\alpha G + w)\Bigr)~.
\end{equation}
In Section~\ref{linearization}, we show that $(\cL -\alpha
\Lambda)^{-1}\cM$ is a bounded operator in $Y$ whose norm is uniformly 
bounded for all $\alpha \in \real$. In Section~\ref{large}, we prove 
that the function $w_\alpha \in Y$ defined by
\begin{equation}\label{walphadef}
  w_\alpha \,=\, -\alpha(\cL - \alpha\Lambda)^{-1}\cM G
\end{equation}
is uniformly bounded in $Y$ for all $\alpha \in \real$, and converges
to some limit $w_\infty$ as $|\alpha| \to \infty$. After these
preliminaries, a standard contraction argument allows to prove that
\reff{weq2} has a unique solution $w^{\alpha,\lambda}$ (in an
appropriate ball in $Y$) if $\lambda$ is sufficiently small, and that
$w^{\alpha,\lambda} = \lambda w_\alpha + \cO(\lambda^2)$.  This is
done in Section~4, where additional properties of the asymmetric
Burgers vortex $\omega^{\alpha,\lambda} = \alpha G +
w^{\alpha,\lambda}$ are also established. Theorem~\ref{thm2} is proved
in Section~\ref{stability} by an energy estimate, using the
observation that the linearization of \reff{omeq2} at the vortex
$\omega^{\alpha,\lambda}$ is a small perturbation of the linear
equation $\partial_t \omega = (\cL - \alpha\Lambda)\omega$ if
$\lambda$ is small.


\section{Linearization in the symmetric case}
\label{linearization}

If $\cL, \Lambda$ are the linear operators defined in 
\reff{LMdef}, \reff{Lambdadef}, we know from \cite{gallay:2004}
that $\cL - \alpha\Lambda$ is invertible in $X$ (with bounded
inverse) for all $\alpha \in \real$. In this section we use the
methods of \cite{gallay:2004} to establish the following result:

\begin{proposition}\label{linbounds}
There exist positive constants $K_1, K_2$ such that, for all 
$\alpha \in \real$, the following inequalities hold:
\begin{eqnarray}\label{XYbd}
  \|(\cL - \alpha\Lambda)^{-1}w\|_Y &\le& K_1 \|w\|_X~, \quad
  \hbox{for all } w \in X~, \\ \label{YYbd}
  \|(\cL - \alpha\Lambda)^{-1}\cM w\|_Y &\le& K_2 \|w\|_Y~, \quad
  \hbox{for all } w \in Y~.
\end{eqnarray}
\end{proposition}

\proof The properties of the linear operator $\cL$ acting on $X$
are easy to establish, because this operator is conjugated to the
Hamiltonian of the harmonic oscillator in $\real^2$. Let $H$ be
the closed subspace of $L^2(\real^2)$ defined by $H = \{f \in 
L^2(\real^2) \,|\, \inttwo G^{1/2}f \d x = 0\}$. 
Consider the linear operator $L : D(L) \to H$ defined by
\begin{eqnarray*}
  D(L) &=& \{f \in H \,|\, \Delta f \in L^2(\real^2)\,,~
  |x|^2 f \in L^2(\real^2)\}~,\\
  L &=& G^{-1/2} (-\cL) G^{1/2} \,=\, -\Delta + \frac{|x|^2}{16}
  - \frac12~.
\end{eqnarray*}
As is well-known, $L$ is self-adjoint in $H$ with spectrum 
$\sigma(L) = \{n/2 \,|\, n = 1,2,3,\dots\}$. In particular, 
$L \ge 1/2$. It follows that there exists $C_1 > 0$ such that, 
for all $f \in H$, 
\begin{equation}\label{bdf1}
  \|L^{-1/2}f\|_{L^2} + \||x| L^{-1/2}f\|_{L^2} + 
  \|\nabla L^{-1/2}f\|_{L^2} \,\le\, C_1 \|f\|_{L^2}~.
\end{equation}
Indeed, setting $g = L^{-1/2}f$, we have
\begin{eqnarray}\nonumber
  \|f\|_{L^2}^2 &=& \|L^{1/2}g\|_{L^2}^2 \,=\, (g,Lg)_{L^2}\\ \label{Lbdd}
  &=& \frac14 \inttwo \Bigl(|\nabla g|^2 + \frac{|x|^2}{16}g^2
  - \frac12 g^2\Bigr)\d x + \frac34 (g,Lg)_{L^2}\\ \nonumber
  &\ge& \frac14 \inttwo \Bigl(|\nabla g|^2 + \frac{|x|^2}{16}g^2
  + g^2\Bigr)\d x~,
\end{eqnarray}
where we have used the fact that $(g,Lg)_{L^2} \ge \half \|g\|_{L^2}^2$.
The bounds \reff{bdf1} also show that $L^{-1/2}$ is a compact operator
in $H$. 

Another useful estimate can be obtained from \reff{bdf1} by a duality
argument. For $j = 1$ or $2$, let $H_j = \{f \in H \,|\, \inttwo 
x_j G^{1/2} f \d x = 0\}$. Then $L^{-1/2}x_j$ extends to a bounded 
operator from $H_j$ into $H$, and for all $f \in H_j$ we have
\begin{equation}\label{bdf2}
  \|L^{-1/2}x_j f\|_{L^2} \,\le\, C_1 \|f\|_{L^2}~.
\end{equation}
Indeed, by density, it is sufficient to prove \reff{bdf2} for 
$f \in H_j \cap \cS(\real^2)$, where $\cS$ denotes the
Schwartz space of test functions. In that case $x_j f \in H$
and for all $\phi \in H$ we have
$$
  |(\phi,L^{-1/2}x_j f)_{L^2}| \,=\, |(x_j L^{-1/2}\phi,f)_{L^2}|
  \,\le\, \|x_j L^{-1/2}\phi\|_{L^2} \|f\|_{L^2} \,\le\, 
  C_1 \|\phi\|_{L^2}\|f\|_{L^2}~,
$$
which proves \reff{bdf2}. 

We now return to the operator $\cL : D(\cL) \to X$ defined by $D(\cL)
= G^{1/2}D(L)$, $X = G^{1/2}H$, and $-\cL \,=\, G^{1/2}L G^{-1/2}$. 
By construction, $\cL$ is selfadjoint in $X$, $-\cL \,\ge\, 1/2$, 
and $(-\cL)^{-1/2}$ is a compact operator in $X$. Using \reff{bdf1}, 
\reff{bdf2} we easily obtain the following additional properties:

\begin{lemma}\label{lemcL}~\\
i) $(-\cL)^{-1/2}$ is a bounded operator from $X$ into $Y$;\\
ii) $(-\cL)^{-1/2}\cM$ extends to a bounded operator from $Y$
into $X$. 
\end{lemma}

\noindent{\bf Proof of Lemma~\ref{lemcL}.} \\
i) Let $w \in X$. Since $(-\cL)^{-1/2} = G^{1/2}L^{-1/2}G^{-1/2}$, 
we have
\begin{eqnarray*}
   \|(-\cL)^{-1/2}w\|_X &=& \|G^{-1/2}(-\cL)^{-1/2}w\|_{L^2}
   \,=\, \|L^{-1/2}G^{-1/2}w\|_{L^2} \\
   &\le& C_1 \|G^{-1/2}w\|_{L^2} \,=\, C_1 \|w\|_X~.
\end{eqnarray*}
Similarly, since $\nabla G^{1/2} = -\frac{x}{4}G^{1/2}$, we find
\begin{eqnarray*}
  \|\nabla (-\cL)^{-1/2}w\|_X &=& \|G^{-1/2}\nabla G^{1/2}
  L^{-1/2} G^{-1/2}w\|_{L^2} \\
  &\le& \|\nabla L^{-1/2} G^{-1/2}w\|_{L^2} + \frac14
  \||x|L^{-1/2} G^{-1/2}w\|_{L^2} \\
  &\le& C_1 \|G^{-1/2}w\|_{L^2} \,=\, C_1 \|w\|_X~.
\end{eqnarray*}
ii) For any $w \in Y$, we have $\|(-\cL)^{-1/2}\cM w\|_X = \half 
\|L^{-1/2} G^{-1/2}(x_1\partial_1{-}x_2\partial_2)w\|_{L^2}$. 
Now $G^{-1/2}\partial_j w \in H_j$ for $j = 1,2$, hence by 
\reff{bdf2}
$$
  \|L^{-1/2}x_j G^{-1/2}\partial_j w\|_{L^2} \,\le\, C_1 
  \|G^{-1/2}\partial_j w\|_{L^2} \,\le\, C_1 \|w\|_Y~.
$$
We conclude that $\|(-\cL)^{-1/2}\cM w\|_X \le C_1\|w\|_Y$. \QED

\medskip
Finally, we consider the operator $\Sigma = (-\cL)^{-1/2}\Lambda
(-\cL)^{-1/2}$, where $\Lambda$ is defined by \reff{Lambdadef}. 
The following properties of $\Sigma$ will be useful:

\begin{lemma}\label{lemSigma}~\\
i) The operator $\Sigma$ is compact in $X$;\\
ii) The operator $\Sigma$ is skew-symmetric in $X$;\\
iii) For any $\alpha \in \real$, the operator $\oone + \alpha\Sigma$
is invertible in $X$ and
\begin{equation}\label{Sigmabd}
  \|(\oone + \alpha\Sigma)^{-1}w\|_X \,\le\, \|w\|_X~, \quad
  \hbox{for all } w \in X~.
\end{equation}
\end{lemma}

\noindent{\bf Proof of Lemma~\ref{lemSigma}.} \\
i) Since $(-\cL)^{-1/2}$ is compact in $X$ and bounded from $X$ 
into $Y$, it suffices to show that $\Lambda : Y \to X$ is bounded. 
If $w \in Y$, then
$$
  \|\vv^G \cdot \nabla w\|_X \,\le\, \|\vv^G\|_{L^\infty}
  \|\nabla w\|_X \,\le\, C \|w\|_Y~.
$$
On the other hand, by H\"older's inequality, we have for all 
$p \in [1,2]$:
\begin{equation}\label{hoelder}
  \|w\|_{L^p} \,=\, \|G^{1/2} G^{-1/2}w\|_{L^p} \,\le\, 
  \|G^{1/2}\|_{L^\frac{2p}{2-p}} \|G^{-1/2}w\|_{L^2} \,=\, C \|w\|_X~.
\end{equation}
If $\vv$ denotes the velocity field obtained from $w$ via the
Biot-Savart law \reff{BS}, the Hardy-Littlewood-Sobolev 
inequality \cite{lieb:1997} implies that $\vv \in L^q(\real^2)$ for 
all $q \in (2,\infty)$, and the following bounds hold:
\begin{equation}\label{HLS}
   \|\vv\|_{L^q} \,\le\, C_p \|w\|_{L^p}~, \quad \hbox{where} 
   \quad 1 < p < 2 \quad \hbox{and} \quad \frac 1q = \frac1p - 
   \frac12~.
\end{equation}
Thus
$$
  \|\vv \cdot \nabla G\|_X \,=\, \|G^{-1/2}\vv \cdot \nabla G\|_{L^2} 
  \,\le\, \|\vv\|_{L^4} \|G^{-1/2} \nabla G\|_{L^4} \,\le\, 
  C \|w\|_{L^{4/3}} \,\le\, C \|w\|_X~.
$$
We conclude that $\|\Lambda w\|_X \le C \|w\|_Y$. 

\noindent ii) It is shown in (\cite{gallay:2004}, Lemma~4.8) that
$(w_1,\Lambda w_2)_X + (\Lambda w_1,w_2)_X = 0$ for all $w_1,w_2 
\in Y$.  Since $(-\cL)^{-1/2}$ is symmetric in $X$, it follows that
$(w_1,\Sigma w_2)_X + (\Sigma w_1,w_2)_X = 0$ for all $w_1,w_2 \in X$.

\noindent iii) By the analytic Fredholm theorem \cite{reed:1980}, we
know that $\oone + \alpha\Sigma$ is invertible in $X$ for all 
$\alpha \in \real$ except perhaps on a discrete set (with no limit
point) where the meromorphic map $\alpha \mapsto (\oone + 
\alpha\Sigma)^{-1}$ has poles. But whenever $\oone + \alpha\Sigma$ is 
invertible we have by ii):
$$
  \|w\|_X^2 \,=\, (w,(\oone{+}\alpha\Sigma)w)_X \,\le\, 
  \|w\|_X \|(\oone{+}\alpha\Sigma)w\|_X \quad \hbox{for all }
  w \in X~,
$$
hence $\|(\oone + \alpha\Sigma)^{-1}w\|_X \le \|w\|_X$. This implies
that $\oone + \alpha\Sigma$ is invertible for all $\alpha \in \real$ 
and that \reff{Sigmabd} holds. \QED

\medskip Equipped with these lemmas, it is now straightforward to
conclude the proof of Proposition~\ref{linbounds}. For any 
$\alpha \in \real$, the formula
$$
  (\cL - \alpha\Lambda)^{-1} \,=\, -(-\cL)^{-1/2} 
  (\oone + \alpha\Sigma)^{-1} (-\cL)^{-1/2} 
$$
shows that $\cL - \alpha\Lambda$ is invertible in $X$ for all
$\alpha \in \real$. The bounds \reff{XYbd} and \reff{YYbd} are
then direct consequences of this identity and Lemmas~\ref{lemcL},
\ref{lemSigma}. \QED


\section{Large Reynolds number asymptotics}
\label{large}

The main goal of this section is to prove that the function $w_\alpha$
defined by \reff{walphadef} is uniformly bounded in the space $Y$ for
all $\alpha \in \real$. From \reff{walphadef} we expect that $w_\alpha
\to \Lambda^{-1}\cM G$ as $|\alpha| \to \infty$, but it is not clear
a priori that this limit makes sense because $\Lambda$ is not an 
invertible operator. Our first result shows that $\cM G$ is indeed
in the range of $\Lambda$ (we recall that $\cS(\real^2)$ is the 
Schwartz space of test functions): 

\begin{proposition}\label{wexist}
There exists $w_\infty \in \cS(\real^2)$ such that $\Lambda w_\infty = 
\cM G$. 
\end{proposition}

\proof By \reff{GvG}, \reff{LMdef} we have $\cM G = \half(x_1\partial_1
- x_2\partial_2)G = -\frac14(x_1^2-x_2^2)G$. Using polar coordinates 
in $\real^2$, we thus find
\begin{equation}\label{MG}
  (\cM G)(r\cos\theta,r\sin\theta) \,=\, -\frac14 r^2 g(r)
  \cos(2\theta)~, 
\end{equation}
where $g(r) = G(r\cos\theta,r\sin\theta) = (4\pi)^{-1}e^{-r^2/4}$. 
As was observed in \cite{robinson:1984,moffatt:1994,prochazka:1995,
gallay:2004}, the operator $\Lambda$ is invariant under rotations 
in the plane, and is therefore block-diagonal in the Fourier basis
$\{e^{in\theta}\}_{n \in \allint}$. We make the following Ansatz:
\begin{equation}\label{ansatz}
  w_\infty(r\cos\theta,r\sin\theta) \,=\,\omega(r) \sin(2\theta)~,
\end{equation}
where $\omega : \real_+ \to \real$ has to be determined. The velocity 
field associated to $w_\infty$ reads
$$
  \vv_\infty \,=\, \frac{2}{r} \Omega(r) \cos(2\theta) \ee_r
  - \Omega'(r) \sin(2\theta)\ee_\theta~,
$$
where $\ee_r$ is the unit vector in the radial direction, $\ee_\theta 
= \ee_r^\perp$, and where $\Omega : \real_+ \to \real$ is the solution 
of the ordinary differential equation
\begin{equation}\label{Omeq1}
  -\frac1r (r\Omega'(r))' + \frac{4}{r^2}\Omega(r) \,=\, \omega(r)~, 
  \quad r > 0~,
\end{equation}
which satisfies the boundary conditions $\Omega(0) = \Omega(+\infty)
= 0$. Using these expressions, we find
\begin{eqnarray}\nonumber
  (\Lambda w_\infty)(r\cos\theta,r\sin\theta) &=& 
  (\vv^G \cdot \nabla w_\infty + \vv_\infty \cdot \nabla G)
  (r\cos\theta,r\sin\theta) \\ \label{Lambdaw} 
  &=& \cos(2\theta) (2\phi(r)\omega(r) - g(r)\Omega(r))~,
\end{eqnarray}
where $\phi(r) = (2\pi r^2)^{-1}(1 - e^{-r^2/4})$. If we compare 
\reff{MG} and \reff{Lambdaw} we obtain the solution
\begin{equation}\label{omOm}
  \omega(r) \,=\, h(r) \Bigl(\Omega(r) - \frac{r^2}4\Bigr)~,
\end{equation}
where $h = g/(2\phi)$, i.e. $h(r) = (r^2/4)(e^{r^2/4}-1)^{-1}$.
Inserting \reff{omOm} into \reff{Omeq1}, we see that $\Omega$ should
satisfy the ordinary differential equation
\begin{equation}\label{Omeq2}
  -\frac1r (r\Omega'(r))' + \Bigl(\frac{4}{r^2} - h(r)\Bigr) 
  \Omega(r) \,=\, -\frac{r^2 h(r)}{4}~, \quad r > 0~,
\end{equation}
together with the boundary conditions $\Omega(0) = \Omega(+\infty)
= 0$. 

\begin{remark}
Equation \reff{Omeq2} was derived and studied numerically by 
Moffatt, Kida and Ohkitani, see Eq.(2.25) in \cite{moffatt:1994}.
The notation used in \cite{moffatt:1994} is $f = \Omega$ and 
$\eta = -h$. 
\end{remark}

To solve \reff{Omeq2}, we first consider the associated homogeneous
equation
\begin{equation}\label{Omeq3}
  -\frac1r (r\Omega'(r))' + \Bigl(\frac{4}{r^2} - h(r)\Bigr) 
  \Omega(r) \,=\, 0~, \quad r > 0~.
\end{equation}
Setting $r = e^{\pm t}$ and $\Omega(r) = F(\pm \log(r))$, this equation
is transformed into
\begin{equation}\label{Feq}
  -F''(t) + (4 - H_\pm(t)) F(t) \,=\, 0~, \quad t \in \real~,
\end{equation}
where $H_\pm(t) = e^{\pm 2t} h(e^{\pm t})$. In particular, $H_\pm(t)$
decays rapidly to zero as $t \to +\infty$. Applying Theorem~3.8.1
in \cite{coddington:1955}, we deduce that \reff{Feq} has a unique
solution $F_\pm(t)$ such that
$$
  \lim_{t \to +\infty} e^{2t} \pmatrix{F_\pm(t) \cr F_\pm'(t)} 
  \,=\, \pmatrix{1 \cr -2}~.
$$

We now define $\psi_+(r) = F_+(\log(r))$ and $\psi_-(r) = 
F_-(-\log(r))$. By construction, $\psi_+, \psi_-$ are the only
solutions of \reff{Omeq3} such that
$$
  \psi_+(r) \,\sim\, \frac{1}{r^2} \quad \hbox{as } r \to +\infty~,
  \qquad \psi_-(r) \,\sim\, r^2 \quad \hbox{as } r \to 0~.
$$
We observe that the ``potential'' $4/r^2   - h(r)$ in
\reff{Omeq3} is strictly positive, because
$$
  \inf_{r > 0} \Bigl(\frac{4}{r^2} - h(r)\Bigr) \,=\, 
  \inf_{z > 0} \Bigl(\frac{1}{z} - \frac{z}{e^z-1}\Bigr) > 0~. 
$$
By the Maximum Principle \cite{protter:1967}, it follows that
$\psi_-'(r) > 0$ and $\psi_+'(r) < 0$ for all $r > 0$. In particular, 
$\psi_+$ and $\psi_-$ are linearly independent, hence there exists
$w_0 > 0$ such that
$$
  W(r) \,=\, \psi_+(r)\psi_-'(r) - \psi_+'(r)\psi_-(r) \,=\, 
  \frac{w_0}{r}~, \quad r > 0~.
$$
Moreover, 
$$
  \psi_+(r) \,\sim\, \frac{w_0}{4r^2} \quad \hbox{as } r \to 0~,
  \qquad \psi_-(r) \,\sim\, \frac{w_0 r^2}{4} \quad \hbox{as } 
  r \to +\infty~.
$$

Using these notations and the ``variation of constants'' formula, 
we obtain the following expression for the solution of \reff{Omeq2}:
\begin{equation}\label{Omsol}
  \Omega(r) \,=\, -\psi_+(r)\int_0^r \frac{z^3}{4 w_0}\psi_-(z)
  h(z)\d z -\psi_-(r)\int_r^\infty \frac{z^3}{4 w_0}\psi_+(z)
  h(z)\d z~, \quad r > 0~.
\end{equation}
It is clear that $\Omega : \real_+ \to \real$ is a smooth function
satisfying 
$$
  \Omega(r) \,\sim\, \left\{ 
  \begin{array}{lcl} \Omega_+ r^2 & as & r \to 0~, \\
  \Omega_- r^{-2} & as & r \to +\infty~, \end{array}\right. 
  \qquad \hbox{where} \quad \Omega_\pm \,=\, -\int_0^\infty 
  \frac{z^3}{4 w_0}\psi_\pm(z)h(z)\d z~.
$$ 
(The values $\Omega_+ \approx -0.38$ and $\Omega_- \approx -17.5$ 
were found numerically in \cite{moffatt:1994}.) Similar estimates
hold for all derivatives. Going back to \reff{omOm}, we see that
$\omega : \real_+ \to \real$ is smooth and rapidly decreasing at
infinity. Moreover $\omega(0) = 0$, $\omega(r) < 0$ for all $r > 0$, 
and it is easy to verify that the Taylor expansion of $\omega(r)$ at 
$r = 0$ contains {\em even} powers of $r$ only. Thus the 
function the function $w_\infty : \real^2 \to \real$ defined by 
\reff{ansatz} is smooth, rapidly decreasing at infinity, and 
satisfies $\Lambda w_\infty = \cM G$ by construction.
\QED

\medskip 
As a consequence we can prove that the function $w_\alpha$ defined
in \reff{walphadef} is uniformly bounded in $Y$:

\begin{corollary}\label{wbound}
There exists $K_3 > 0$ such that $\|w_\alpha\|_Y \le K_3|\alpha|/(1{+}
|\alpha|)$ for all $\alpha \in \real$. 
\end{corollary}

\proof From \reff{walphadef} and \reff{YYbd} we know that 
$\|w_\alpha\|_Y \le K_2 |\alpha| \|G\|_Y$. On the other hand, in 
view of Proposition~\ref{wexist}, we have $(\cL - \alpha\Lambda) 
w_\alpha = -\alpha \cM G = (\cL - \alpha\Lambda)w_\infty - \cL w_\infty$, 
hence
\begin{equation}\label{walphaexp}
  w_\alpha \,=\, w_\infty - (\cL - \alpha\Lambda)^{-1} \cL w_\infty~.
\end{equation}
Using \reff{XYbd}, we infer that $\|w_\alpha\|_Y \le \|w_\infty\|_Y 
+ K_1 \|\cL w_\infty\|_X$. Combining both estimates we obtain the 
desired result. \QED

\medskip
A more detailed analysis reveals that $w_\infty$ is indeed the limit
of $w_\alpha$ as $|\alpha| \to \infty$:

\begin{proposition}\label{wlimit}
There exists $K_4 > 0$ such that $\|w_\alpha - w_\infty\|_Y \le 
K_4/(1{+}|\alpha|)$ for all $\alpha \in \real$. 
\end{proposition}

\proof By \reff{walphaexp} and \reff{XYbd}, we have $\|w_\alpha - 
w_\infty\|_Y \le K_1 \|\cL w_\infty\|_X$. On the other hand, 
proceeding exactly as in the proof of Proposition~\ref{wexist},
it is straightforward to show that there exists $z_\infty \in 
\cS(\real^2)$ such that $\Lambda z_\infty = \cL w_\infty$ (the details
are left to the reader). Then, using \reff{walphaexp}, we find for 
all $\alpha \neq 0$
$$
  w_\alpha \,=\, w_\infty + \frac{1}\alpha z_\infty - \frac{1}\alpha
  (\cL - \alpha\Lambda)^{-1} \cL z_\infty~,
$$
hence $\|w_\alpha - w_\infty\|_Y \le |\alpha|^{-1}(\|z_\infty\|_Y +
K_1 \|\cL z_\infty\|_X)$. Combining both estimates we obtain the 
desired result. \QED


\section{Existence of asymmetric vortices}
\label{existence} 
 
As is explained in the introduction, we shall prove the existence of 
a stationary solution of \reff{omeq2} by solving Eq.\reff{weq2}, 
namely
\begin{equation}\label{weq3}  
  w \,=\, \lambda w_\alpha + (\cL - \alpha\Lambda)^{-1}
  \Bigl(\vv\cdot\nabla w - \lambda\cM w \Bigr)~,
\end{equation}
where $w_\alpha$ is defined in \reff{walphadef} and (as usual) $\vv$ 
denotes the velocity field obtained from $w$ by the Biot-Savart law
\reff{BS}. To bound the nonlinear term in \reff{weq3}, we use the
following bilinear estimate:

\begin{lemma}\label{bilinear}  
There exists $K_5 > 0$ such that, if $w, \tilde w \in Y$ and if $\vv$ 
is the velocity field obtained from $w$ by the Biot-Savart law, then
$\vv \cdot \nabla \tilde w \in X$ and
$$
  \|\vv \cdot \nabla \tilde w\|_X \,\le\, K_5 \|w\|_Y 
  \|\tilde w\|_Y~.
$$
\end{lemma}

\proof We first observe that $w \in L^p(\real^2)$ for all $p \in 
[1,+\infty)$, and that $\|w\|_{L^p} \le C_p \|w\|_Y$. Indeed this is
proved in \reff{hoelder} if $p \in [1,2]$, and if $p \ge 2$ this 
follows from the embeddings $Y \hookrightarrow H^1(\real^2) 
\hookrightarrow L^p(\real^2)$. Next, by the Calder\'on-Zygmund 
inequality \cite{stein:1970}, the velocity field $\vv$ associated 
to $w$ satisfies $\|\nabla \vv\|_{L^p} \le C_p \|w\|_{L^p}$ for 
$1 < p < \infty$. Using the Gagliardo-Nirenberg inequality 
\cite{nirenberg:1959}, we deduce that $\vv \in L^\infty(\real^2)$ 
and that
$$
  \|\vv\|_{L^\infty} \,\le\, C \|\nabla \vv\|_{L^3}^{1/2} 
  \|\vv\|_{L^6}^{1/2} \,\le\, C \|w\|_{L^3}^{1/2}\|w\|_{L^{3/2}}^{1/2}
  \,\le\, C \|w\|_Y~,
$$
where we also used the Hardy-Littlewood-Sobolev inequality \reff{HLS}
with $q = 6$, $p = 3/2$. We conclude that $\vv \cdot \nabla \tilde w
\in X$ and
$$
  \|\vv \cdot \nabla \tilde w\|_X \,\le\, \|\vv\|_{L^\infty}
  \|\nabla \tilde w\|_X \,\le\, C \|w\|_Y \|\tilde w\|_Y~,
$$
which is the desired result. \QED

\medskip 
We can now prove the main result of this section:

\begin{proposition}\label{fixedpoint}
Choose $\lambda_0 \in (0,1)$ such that
\begin{equation}\label{lambda0}
  2 K_2 \lambda_0 \,<\, 1~, \quad \hbox{and}\quad 
  16 K_1 K_3 K_5 \lambda_0 \,\le\, 1~.
\end{equation}
Then, for all $\lambda \in [0,\lambda_0]$ and all $\alpha \in \real$, 
equation~\reff{weq3} has a unique solution $w^{\alpha,\lambda} \in Y$
such that $\|w^{\alpha,\lambda}\|_Y \le (4K_1 K_5)^{-1}$. This solution
depends smoothly on the parameters $(\alpha,\lambda) \in \real \times
[0,\lambda_0]$ and satisfies $\|w^{\alpha,\lambda}\|_Y \le 4 K_3
\lambda |\alpha|/(1{+}|\alpha|)$.
\end{proposition}

\proof Fix $\lambda \in [0,\lambda_0]$, $\alpha \in \real$, and choose
$r > 0$ such that
\begin{equation}\label{rdef}
  4 K_3^\alpha \lambda \,\le\, r \,\le\, \frac{1}{4K_1 K_5}~, \quad
  \hbox{where} \quad K_3^\alpha \,=\, \frac{K_3|\alpha|}{1+|\alpha|}~.
\end{equation}
Let $F^{\alpha,\lambda} : Y \to Y$ be the quadratic map defined by 
$$
  F^{\alpha,\lambda}(w) \,=\, \lambda w_\alpha + (\cL - \alpha\Lambda)^{-1}
  \Bigl(\vv\cdot\nabla w - \lambda\cM w \Bigr)~, \quad w \in Y~.
$$
Using Proposition~\ref{linbounds}, Corollary~\ref{wbound} and 
Lemma~\ref{bilinear}, we find for all $w \in Y$:
\begin{equation}\label{Fmap}
  \|F^{\alpha,\lambda}(w)\|_Y \,\le\, K_3^\alpha \lambda + K_1 K_5 
  \|w\|_Y^2 + K_2\lambda \|w\|_Y~.
\end{equation}
Similarly, for all $w, \tilde w \in Y$:
\begin{equation}\label{Fcontract}
  \|F^{\alpha,\lambda}(w) - F^{\alpha,\lambda}(\tilde w)\|_Y \,\le\, 
  \Bigl(K_1 K_5 (\|w\|_Y + \|\tilde w\|_Y) +  K_2\lambda\Bigr) 
  \|w - \tilde w\|_Y~.
\end{equation}
Let $B_r = \{w \in Y\,|\, \|w\|_Y \le r\}$. It follows from
\reff{lambda0}, \reff{rdef}, and \reff{Fmap} that $F^{\alpha,\lambda}$ maps 
$B_r$ into itself, because $K_3^\alpha \lambda + K_1 K_5 r^2 + 
K_2 \lambda r \le r/4 + r/4 + r/2 = r$. Similarly \reff{Fcontract}
implies for all $w, \tilde w \in B_r$:
$$
  \|F^{\alpha,\lambda}(w) - F^{\alpha,\lambda}(\tilde w)\|_Y \,\le\, 
  \kappa \|w - \tilde w\|_Y~, \quad \hbox{where} \quad \kappa \,=\,
  \half + K_2 \lambda_0 < 1~.
$$
By the contraction mapping theorem, $F^{\alpha,\lambda}$ has thus a
unique fixed point in $B_r$, which we denote by $w^{\alpha,\lambda}$.
It remains to show that $w^{\alpha,\lambda}$ is a smooth function of
$(\alpha,\lambda) \in \real \times [0,\lambda_0]$. But this is a
direct consequence of the implicit function theorem, because the map
$F^{\alpha,\lambda}$ depends smoothly on $(\alpha,\lambda)$ and the
differential
$$
  D_w F^{\alpha,\lambda}(w) \,=\, \tilde w \mapsto 
  (\cL - \alpha\Lambda)^{-1} \Bigl(\tilde \vv\cdot\nabla w +
  \vv\cdot\nabla \tilde w - \lambda\cM \tilde w \Bigr)
$$
satisfies $\|D_w F^{\alpha,\lambda}(w)\| \le \kappa$ for all $w \in
B_r$ and all $(\alpha,\lambda) \in \real \times [0,\lambda_0]$.
Thus $\oone - D_w F^{\alpha,\lambda}(w)$ is invertible at $w = 
w^{\alpha,\lambda}$, and the desired conclusion follows from the 
implicit function theorem. \QED

\medskip
Theorem~\ref{thm1} is an immediate consequence of 
Proposition~\ref{fixedpoint}: we just set $\omega^{\alpha,\lambda} = 
\alpha G + w^{\alpha,\lambda}$, $\uu^{\alpha,\lambda} = \alpha \vv^G +
\vv^{\alpha,\lambda}$ (where $\vv^{\alpha,\lambda}$ is the velocity 
field obtained from $w^{\alpha,\lambda}$ by the Biot-Savart law), 
and $K_0 \,=\, (4K_1 K_5)^{-1}$. By construction, $\omega^{\alpha,
\lambda}$ is a stationary solution of \reff{omeq2} satisfying 
the conclusions of the theorem. \QED

In the rest of this section, we establish a few additional properties
of the asymmetric vortex $\omega^{\alpha,\lambda}$:

\medskip\noindent{\bf 1)} {\em Expansion in $\lambda$.} There exists
$K_6 > 0$ such that, for all $(\alpha,\lambda) \in \real \times 
[0,\lambda_0]$:
\begin{equation}\label{lambdaexp}
  \|\omega^{\alpha,\lambda} - \alpha G -\lambda w_\alpha\|_Y \,\le\, 
  \frac{K_6 |\alpha|\lambda^2}{1+|\alpha|}~.
\end{equation}
Indeed, using the notations of Proposition~\ref{fixedpoint}, we have
$w^{\alpha,\lambda} \in B_r$ with $r = 4 K_3^\alpha \lambda$. 
As $w^{\alpha,\lambda}$ is a solution of \reff{weq3}, we obtain
$$
  \|w^{\alpha,\lambda} -\lambda w_\alpha\|_Y \,\le\, K_1 K_5 r^2 
  + K_2 \lambda r \,\le\, \frac{K_6 |\alpha|\lambda^2}{1+|\alpha|}~.
$$

\medskip\noindent{\bf 2)} {\em Large Reynolds number asymptotics.}
Combining \reff{lambdaexp} and Proposition~\ref{wlimit}, we find
for $|\alpha| \ge 1$:
\begin{equation}\label{largeR}
  \Bigl\|\frac{1}{\alpha}\omega^{\alpha,\lambda} - G -\frac{\lambda}{\alpha}
  w_\infty\Bigr\|_Y \,\le\, \frac{K_6 \lambda^2}{|\alpha|} + 
  \frac{K_4 \lambda}{|\alpha|^2}~. 
\end{equation}
In agreement with \cite{moffatt:1994}, we see that the leading
correction to the Gaussian profile $G$ is $(\lambda/\alpha)w_\infty$, 
and that the higher order corrections are proportional to $\lambda^2/R$
and $\lambda/R^2$, where $R = |\alpha|$ is the Reynolds number. 

\medskip\noindent{\bf 3)} {\em Small Reynolds number asymptotics.}
Since $\cL\cM G = \cL^{-1}\cM G = -\cM G$, it follows from
\reff{walphadef} that $w_\alpha = \alpha \cM G + \alpha^2 
(\cL - \alpha\Lambda)^{-1}\Lambda\cM G$. Replacing into 
\reff{lambdaexp} and using Proposition~\ref{linbounds}, we obtain
for $|\alpha| \le 1$:
\begin{equation}\label{smallR1}
  \|\omega^{\alpha,\lambda} - \alpha G -\lambda \alpha \cM G\|_Y 
  \,\le\, K_6 |\alpha|\lambda^2 + K_1 \alpha^2 \lambda 
  \|\Lambda \cM G\|_X~.
\end{equation}
In fact, if we proceed as in (\cite{gallay:2005}, Section~2), 
this result can be improved as follows: there exists $K_7 > 0$
such that, for $|\alpha| \le 1$,
\begin{equation}\label{smallR2}
  \|\omega^{\alpha,\lambda} - \alpha G_\lambda\|_Y \,\le\, 
  K_7 \alpha^2 \lambda~,
\end{equation}
where
$$
  G_{\lambda}(x) \,=\, \frac{\sqrt{1-\lambda^2}}{4 \pi} 
  \,\exp\Bigl(-\frac{1+\lambda}{4}\,x_1^2 -\frac{1-\lambda}{4}\,x_2^2 
  \Bigr)~,  \quad x \in \real^2~.
$$
Remark that $(\cL + \lambda\cM)G_\lambda = 0$ and $\inttwo G_\lambda
\d x = 1$. Since $G_\lambda = G + \lambda\cM G + \cO(\lambda^2)$, 
we see that \reff{smallR1} is compatible with \reff{smallR2}. 

\medskip\noindent{\bf 4)} {\em Smoothness in $x$.}
Standard elliptic estimates imply that $\omega^{\alpha,\lambda}(x)$
is a smooth function of $x \in \real^2$ for any $(\alpha,\lambda) \in 
\real \times [0,\lambda_0]$. Indeed, since $\omega^{\alpha,\lambda}$
is a stationary solution of \reff{omeq2}, we have
$$
  \omega^{\alpha,\lambda} \,=\, (\cL + \lambda\cM)^{-1} \nabla\cdot
  \left(\uu^{\alpha,\lambda}\omega^{\alpha,\lambda}\right)~.
$$
It is not difficult to prove that the linear operator $(\cL + 
\lambda\cM)^{-1} \nabla$ is regularizing, hence a bootstrap argument
shows that $\omega^{\alpha,\lambda} \in H^k(\real^2)$ for all $k \in
\intplus$. One can also prove that all derivatives decay rapidly
at infinity, so that $\omega^{\alpha,\lambda} \in \cS(\real^2)$. 

\medskip\noindent{\bf 5)} {\em Positivity.} It follows from the
parabolic Maximum Principle \cite{protter:1967} that
$\omega^{\alpha,\lambda}(x) > 0$ for all $x \in \real^2$ if $\alpha >
0$. (Similarly, $\omega^{\alpha,\lambda}(x) < 0$ if $\alpha < 0$, and
we already know from \reff{smallR1} that $\omega^{\alpha,\lambda}
\equiv 0$ if $\alpha = 0$.) Indeed, arguing as in the symmetric
case $\lambda = 0$, it is not difficult to show that, for any initial
data $\omega_0 \in L^1(\real^2) \cap C^0(\real^2)$, Eq.\reff{omeq2}
has a unique global solution $\omega \in C^0([0,+\infty),L^1(\real^2))$
which satisfies $\inttwo \omega(x,t)\d x = \inttwo \omega_0\d x$ 
for all $t \ge 0$. If $\omega_0 \ge 0$ and $\omega_0$ is not
identically zero, the Maximum Principle implies that $\omega(x,t) > 0$
for all $x \in \real^2$ and all $t > 0$, see e.g. (\cite{gallay:2004}, 
Section~2.3). Conversely, if $\omega_0$ has non-constant sign, the
norm $\|\omega(\cdot,t)\|_{L^1}$ is strictly decreasing in
time (for $t > 0$ sufficiently small), see (\cite{gallay:2004}, 
Section~3.1). Since $\omega^{\alpha,\lambda} \in L^1(\real^2) \cap 
C^0(\real^2)$ is a stationary solution of \reff{omeq2} satisfying 
$\inttwo \omega^{\alpha,\lambda}\d x = \alpha$, the properties above
imply that $\omega^{\alpha,\lambda}(x) > 0$ for all $x \in \real^2$ if 
$\alpha > 0$. 


\section{Stability of asymmetric vortices}
\label{stability} 
 
In this final section, we show that the asymmetric Burgers vortex
$\omega^{\alpha,\lambda}$ constructed in Section~\ref{existence} is a
stable solution of \reff{omeq2} with respect to perturbations in $X$,
provided $\lambda > 0$ is sufficiently small. Fix $\alpha \in \real$,
$\lambda \in [0,\lambda_0]$, and consider solutions of \reff{omeq2} of
the form $\omega = \omega^{\alpha,\lambda} + \tilde\omega$, $\uu =
\uu^{\alpha,\lambda} + \tilde \uu$, where $\tilde\uu$ is the velocity
field obtained from $\tilde \omega$ by the Biot-Savart law \reff{BS}.
Then $\tilde \omega$ satisfies the equation
\begin{equation}\label{tomeq1}
  \partial_t \tilde\omega + \uu^{\alpha,\lambda} \cdot \nabla 
  \tilde \omega + \tilde\uu \cdot \nabla \omega^{\alpha,\lambda}
  + \tilde\uu \cdot\nabla \tilde\omega \,=\,  
  (\cL + \lambda \cM) \tilde\omega~.
\end{equation}
If we further decompose $\omega^{\alpha,\lambda} = \alpha G + 
w^{\alpha,\lambda}$, $\uu^{\alpha,\lambda} = \alpha \vv^G + 
\vv^{\alpha,\lambda}$, this equation becomes
\begin{equation}\label{tomeq2}
  \partial_t \tilde\omega + \tilde\uu \cdot\nabla \tilde\omega \,=\,  
  \Bigl(\cL + \lambda \cM -\alpha\Lambda -\cN^{\alpha,\lambda}\Bigr) 
  \tilde\omega~,
\end{equation}
where $\Lambda$ is defined in \reff{Lambdadef} and $\cN^{\alpha,\lambda}$
is the integro-differential operator defined by
$$ 
  \cN^{\alpha,\lambda}\tilde \omega \,=\,\vv^{\alpha,\lambda} \cdot 
  \nabla \tilde\omega + \tilde\uu \cdot \nabla w^{\alpha,\lambda}~.
$$
It is easy to show that the Cauchy problem for \reff{tomeq1} or 
\reff{tomeq2} is locally well-posed in the space $X$. In the 
symmetric case $\lambda = 0$, this is proved in \cite{gallay:2002a}
using a larger function space (with polynomial instead of Gaussian
weight), and the same arguments apply here with straightforward 
modifications. Our goal is to control the behavior of the solutions
of \reff{tomeq2} in a neighborhood of the origin. An energy estimate
yields the following result:

\begin{proposition}\label{energy}
There exist positive constants $K_8, K_9$ (independent of 
$\alpha$ and $\lambda$) such that, for any $\delta \in (0,1)$, 
any solution solution of \reff{tomeq2} in $X$ satisfies 
$$
  \frac{\D}{\D t}\|\tilde \omega(t)\|_X^2 \,\le\, 
  - (1{-}\delta) \|\tilde \omega(t)\|_X^2
$$
whenever $K_8\lambda + K_9 \|\tilde \omega(t)\|_X \le \delta$. 
\end{proposition}

\proof From \reff{tomeq2} we obtain
\begin{equation}\label{energyeq}
  \frac12 \frac{\D}{\D t}\|\tilde \omega\|_X^2 \,=\, \Bigl(\tilde
  \omega,(\cL+\lambda\cM-\alpha\Lambda-\cN^{\alpha,\lambda})
  \tilde\omega\Bigr)_{\!X} - (\tilde\omega,\tilde\uu \cdot\nabla
  \tilde\omega)_X~.
\end{equation}
To simplify the subsequent expressions, we define $f = G^{-1/2}\tilde
\omega$, so that $\|\tilde\omega\|_X = \|f\|_{L^2}$, and we introduce
the quadratic form
$$
  E(f) \,=\, \inttwo \Bigl(|\nabla f|^2 + \frac{|x|^2}{16}f^2
  \Bigr)\d x~.
$$
Obviously $\|f\|_{L^2}^2 \le C E(f)$, and we have $\|f\|_{L^4}^2 
\le C \|\nabla f\|_{L^2}^{3/2}\|xf\|_{L^2}^{1/2} \le C E(f)$, 
see \cite{caffarelli:1984}. 

Proceeding as in \reff{Lbdd}, we obtain for any $\delta \in (0,1)$:
$$
  (\tilde\omega,\cL\tilde\omega)_X \,=\, -(f,Lf)_{L^2}
  \,\le\, -\frac{1{-}\delta}{2}\|f\|_{L^2}^2 -\frac{\delta}{2}
  E(f)~.
$$
(The bound \reff{Lbdd} was the particular case $\delta = 1/2$.) 
Next, since $\cM = \half(x_1\partial_1-x_2\partial_2)$ and 
$G^{-1/2}\nabla \tilde\omega = \nabla f - \frac{x}{4}f$, we find
$$
  (\tilde\omega,\cM\tilde\omega)_X \,\le\, \|x\tilde\omega\|_X
  \|\nabla\tilde\omega\|_X \,=\, \|xf\|_{L^2}\,\|\nabla f - 
  \textstyle{\frac{x}{4}}f\|_{L^2} \,\le\, C E(f)~.
$$
Moreover, as was observed in the proof of Lemma~\ref{lemSigma}, 
the operator $\Lambda$ is skew-symmetric in $X$, hence 
$(\tilde\omega,\Lambda\tilde\omega)_X = 0$. 

We now bound the nonlinear term in \reff{energyeq}. Integrating by 
parts and using the relation $\nabla G^{-1} = \frac{x}{2}G^{-1}$,
we obtain
$$
  (\tilde\omega,\tilde\uu\cdot\nabla\tilde\omega)_X \,=\, 
  \inttwo G^{-1}\tilde\omega (\tilde\uu\cdot\nabla\tilde\omega)
  \d x \,=\,-\frac14 \inttwo G^{-1}(x\cdot\tilde\uu)\tilde
  \omega^2 \d x~.
$$
Applying H\"older's inequality, we find
$$
  |(\tilde\omega,\tilde\uu\cdot\nabla\tilde\omega)_X| \,\le\, 
  \frac14 \inttwo |\tilde\uu| |xf| |f| \d x \,\le\, 
  \frac14 \|\tilde\uu\|_{L^4} \|f\|_{L^4} \|xf\|_{L^2} 
  \,\le\, C \|\tilde\omega\|_X E(f)~,
$$
where we have used the bound $\|\tilde\uu\|_{L^4} \le C 
\|\tilde\omega\|_{L^{4/3}} \le C\|\tilde\omega\|_X$, see 
\reff{hoelder} and \reff{HLS}. 

It remains to bound $(\tilde\omega,\cN^{\alpha,\lambda}\tilde\omega)_X
= (\tilde\omega,\vv^{\alpha,\lambda}\cdot\nabla\tilde\omega)_X
+ (\tilde\omega,\tilde\uu\cdot\nabla w^{\alpha,\lambda})_X$. The first
term in this sum can be estimated in the same way as the nonlinear 
term above, namely $|(\tilde\omega,\vv^{\alpha,\lambda}\cdot
\nabla\tilde\omega)_X| \le C \|w^{\alpha,\lambda}\|_X E(f)$. 
For the second term, we argue differently:
\begin{eqnarray*}
  |(\tilde\omega,\tilde\uu\cdot\nabla w^{\alpha,\lambda})_X| &\le& 
  \inttwo G^{-1}|\tilde \omega| |\tilde \uu| |\nabla 
  w^{\alpha,\lambda}|\d x \\ 
  &\le& \|\nabla w^{\alpha,\lambda}\|_X \|\tilde \uu\|_{L^4} 
  \|f\|_{L^4} \,\le\, C \|w^{\alpha,\lambda}\|_Y
  E(f)~,
\end{eqnarray*}
where we have used $\|\tilde\uu\|_{L^4} \le C\|\tilde\omega\|_{L^{4/3}}
\le C\|f\|_{L^2} \le C E(f)^{1/2}$. Thus, using 
Proposition~\ref{fixedpoint} to bound $w^{\alpha,\lambda}$, 
we conclude that $|(\tilde\omega,\cN^{\alpha,\lambda}\tilde\omega)_X| 
\le C \|w^{\alpha,\lambda}\|_Y E(f) \le C\lambda E(f)$.
Summarizing, we have shown that
$$
  \frac12 \frac{\D}{\D t}\|\tilde\omega\|_X^2 \,\le\, 
  -\frac{1{-}\delta}{2}\|\tilde\omega\|_X^2 + 
  \frac12 \Bigl(K_8 \lambda + K_9 \|\tilde\omega\|_X - \delta\Bigr)
  E(G^{-1/2}\tilde\omega)~,
$$
for some constants $K_8,K_9$ independent of $\alpha$, $\lambda$ 
and of the solution $\tilde \omega$ of \reff{tomeq2}. This concludes 
the proof. \QED

\medskip
Theorem~\ref{thm2} is now a direct consequence of 
Proposition~\ref{energy}. Given $\mu \in (0,1/2)$, it suffices to
choose $\lambda_1 \in (0,\lambda_0]$ such that $K_8 \lambda_1 \le
\delta/2$ and $\epsilon > 0$ such that $K_9\epsilon < \delta/2$, 
where $\delta = 1 -2\mu$. If $\tilde \omega_0 \in X$ satisfies
$\|\tilde\omega_0\|_X \le \epsilon$ and if $\tilde \omega \in 
C^0([0,T^*),X)$ denotes the (maximal) solution of \reff{tomeq2}
with initial data $\tilde\omega_0$, we define
$$
  T \,=\, \sup\Bigl\{t \in (0,T^*)\,\Big|\, K_9 \|\tilde\omega(s)\|_X
  \le \delta/2 \hbox{ for }0 \le s \le t\Bigr\} \,\in\, (0,T^*]~.
$$
If $T < +\infty$ then Proposition~\ref{energy} implies that 
$\|\tilde\omega(t)\|_X \le e^{-\mu t}\|\tilde\omega_0\|_X$ for 
$t \in [0,T]$, hence $K_9 \|\tilde\omega(t)\|_X \le K_9\epsilon < 
\delta/2$ for $t \in [0,T]$, which contradicts the definition of $T$.
Thus, we must have $T = T^* = +\infty$. This means that the solution
$\tilde\omega$ is globally defined, and by Proposition~\ref{energy}
$\|\tilde\omega(t)\|_X \le e^{-\mu t}\|\tilde\omega_0\|_X$ for all 
$t \ge 0$. \QED

We conclude with a few remarks on the basin of attraction of the 
asymmetric vortex and the decay rate in time of the perturbations:

\medskip\noindent{\bf a)} Proposition~\ref{energy} shows that the
decay rate in time of pertubations of the asymmetric vortex
$\omega^{\alpha,\lambda}$ is bounded from below by $\mu = \half
(1-K_8\lambda) = \half - \cO(\lambda)$, uniformly in $\alpha \in
\real$. This is consistent with the information we have on the 
spectrum of the linearized operator
$$
  \cL^{\alpha,\lambda} \,=\, \cL + \lambda\cM - \alpha\Lambda 
  -\cN^{\alpha,\lambda}~,
$$
acting on the space $X$. Indeed, differentiating the identity 
$(\cL + \lambda\cM)\omega^{\alpha,\lambda} \,=\, \uu^{\alpha,\lambda}
\cdot\nabla\omega^{\alpha,\lambda}$ with respect to $x_1$ and 
$x_2$ we obtain
$$
  \cL^{\alpha,\lambda}(\partial_1 \omega^{\alpha,\lambda}) \,=\,
  -\frac{1{+}\lambda}{2}(\partial_1 \omega^{\alpha,\lambda})~, \quad
  \cL^{\alpha,\lambda}(\partial_2 \omega^{\alpha,\lambda}) \,=\,
  -\frac{1{-}\lambda}{2}(\partial_2 \omega^{\alpha,\lambda})~.
$$
In particular, since $\partial_2 \omega^{\alpha,\lambda} \in X$, we
see that $-\frac{1{-}\lambda}{2}$ is always an eigenvalue of
$\cL^{\alpha,\lambda}$, hence $\mu \le \frac{1{-}\lambda}{2}$.
Numerical calculations by Prochazka and Pullin \cite{prochazka:1998}
seem to indicate that $-\frac{1{-}\lambda}{2}$ is always the largest
eigenvalue of $\cL^{\alpha,\lambda}$ in $X$ (for any $\lambda$).  If
this was true our arguments could be extended to prove existence and
stability of asymmetric Burgers vortices for all $\lambda \in [0,1)$.
Put another way, our current limitation on the range of the asymmetry
parameter is only due to the fact that we do not know how to 
control the eigenvalues of $\cL^{\alpha,\lambda}$ (except of course
for small $\lambda$). 

\medskip\noindent{\bf b)} A remarkable feature of our stability 
result (Theorem~\ref{thm2}) is that it holds uniformly for all 
$\alpha \in \real$. In particular, this implies a uniform upper
bound on the eigenvalues of the linearized operator 
$\cL^{\alpha,\lambda}$. This is definitely compatible with the 
numerical observations of Prochazka and Pullin \cite{prochazka:1995},
but these calculations suggest that our result is perhaps not 
optimal for large Reynolds numbers. According to \cite{prochazka:1995}
we expect that the eigenvalues that are not frozen by symmetries 
have a real part that converges to $-\infty$ as $|\alpha| \to \infty$,
which could imply a faster decay rate $\mu$ and a larger basin of
attraction $\epsilon$ for large Reynolds numbers. A mathematical 
understanding of these numerical observations is still lacking. 

\medskip\noindent{\bf c)} We chose to consider perturbations 
$\tilde \omega$ in the weighted space $X$ because of the
``miraculous'' fact that the operator $\Lambda$ is skew-symmetric
in that space. This is why Proposition~\ref{energy} holds 
uniformly for all $\alpha \in \real$. However the space $X$ is
relatively small since its elements are forced to decay rapidly
at infinity in space. Extending the methods developed in 
\cite{gallay:2004} for the symmetric case $\lambda = 0$, it is not
difficult to show that the asymmetric Burgers vortices are also
stable with respect to perturbations in weighted $L^2$ spaces 
with {\em polynomial} (instead of Gaussian) weight. The decay rate
in time of the perturbations is still uniform in $\alpha$, but the
size of the basin of attraction is a priori not. 

\medskip\noindent{\bf d)} In the symmetric case $\lambda = 0$ it
is shown in \cite{gallay:2004} that the Burgers vortex $\alpha G$ 
is the unique stationary solution of \reff{omeq2} such that 
$\omega \in L^1(\real^2)$ and $\inttwo \omega\d x = \alpha$. 
Moreover, any solution $\omega \in C^0([0,+\infty),L^1(\real^2))$
of \reff{omeq2} such that $\inttwo \omega(x,t)\d x = \alpha$ 
converges to $\alpha G$ in $L^1(\real^2)$ as $t \to +\infty$. 
We do not know if such global results can be extended to the
nonsymmetric case $\lambda > 0$. 

\bigskip\noindent
{\bf Acknowledgements.}
A part of this work was completed when CEW was a visitor at Institut
Fourier, Universit\'e de Grenoble I, whose hospitality is gratefully
acknowledged. The research of CEW is supported in part by the NSF
through grant DMS-0405724, and the work of ThG is supported by the 
ACI ``Structure and dynamics of nonlinear waves'' of the French 
Ministry of Research. 


\bibliographystyle{plain}

\end{document}